# A note on the supremum of a stable process

R. A. DONEY

The University of Manchester UK

If $X$ is a spectrally positive stable process of index $\alpha \in (1,2)$ whose Lévy measure has density $cx^{-\alpha-1}$ on $(0,\infty)$, and $S_1 = \sup_{0<t\leq 1} X_t$, it is known that $P(S_1 > x) \backsim c\alpha^{-1}x^{-\alpha}$ as $x \to \infty$. It is also known that $S_1$ has a continuous density, $s$ say. The point of this note is to show that $s(x) \backsim cx^{-(\alpha+1)}$ as $x \to \infty$.

## 1. Introduction

In connection with an optimal stopping problem, the following question was posed recently by Robert Dalang [4]. Let $X$ be a strictly stable, spectrally positive process of index $\alpha \in (1,2)$, whose Lévy measure has density $\nu(x) = cx^{-(\alpha+1)}$, $x > 0$. Then if $S$ is its maximum process it is known that both $X_1$ and $S_1$ have continuous density functions, $f$ and $s$ say. From (14.34), p 88 in [7] we know that

$$f(x) \backsim cx^{-(\alpha+1)} \text{ as } x \to \infty, \tag{1}$$

and from Proposition 4, p 221 in [1]

$$P(S_1 > x) \backsim P(X_1 > x) \backsim c\alpha^{-1}x^{-\alpha} \text{ as } x \to \infty. \tag{2}$$

The question is whether in fact we can improve this to

$$s(x) \backsim cx^{-(\alpha+1)} \text{ as } x \to \infty. \tag{3}$$

Of course, if we knew in advance that $s$ was ultimately monotone, this would follow immediately from (2), but as we don't have this information, we have to use a different method.

In [2] it was observed that the fact that in this case the Wiener-Hopf factorisation takes a semi-explicit form can be exploited to give the following information. Without loss of generality take $c = 1/\Gamma(-\alpha)$; (in other cases the results follow by scaling); then

$$E(e^{-\lambda S_1}) = \frac{\alpha}{\Gamma(1/\alpha)} e^{\lambda^\alpha} \int_\lambda^\infty e^{-y^\alpha} dy, \tag{4}$$

$$s(x) = \sum_1^\infty a_n x^{\alpha n - 2}, \text{ with } a_n = \frac{1}{\Gamma(\alpha n - 1)\Gamma(\alpha^{-1} + 1 - n)}, \tag{5}$$

and

$$s(x) = \frac{I_1(x)}{\pi \Gamma(\alpha^{-1})} + \frac{I_2(x)}{\pi}, \quad \text{where} \tag{6}$$

$$I_1(x) = \int_0^\infty \int_0^{t^\alpha} \frac{e^{y\cos(\alpha\pi/2)}}{(t^\alpha - y)^{1-1/\alpha}} \sin(y \sin(\alpha\pi/2) + tx) dy dt,$$

$$I_2(x) = \int_0^\infty e^{t^\alpha \cos(\alpha\pi/2)} \cos(t^\alpha \sin(\alpha\pi/2) + tx) dt.$$

(Note that (5) does NOT follow by expanding the RHS of (4) in negative powers of $\lambda$ :see [2] for a complete proof.)

In this note we show that, despite its complicated form, (6) can be exploited to answer the question in the affirmative. The first point is that, after writing $I_1(x)$ in the form

$$I_1(x) = \int_0^\infty g_1(t) \cos tx \, dt + \int_0^\infty g_2(t) \sin tx \, dt, \tag{7}$$

it is possible to show that both $g_1$ and $g_2$ are three times differentiable and satisfy

$$\int_0^\infty \{|g_1'''(t)| + |g_2'''(t)|\} dt < \infty. \tag{8}$$

Then three integrations by part show that $I_1(x) = O(x^{-3})$ as $x \to \infty$.

Next, we write

$$I_2(x) = \int_0^\infty h_1(t) \cos tx \, dt + \int_0^\infty h_2(t) \sin tx \, dt, \tag{9}$$

and note that, although both $h_1$ and $h_2$ are three times differentiable, their third derivatives have non-integrable singularities at zero. In fact there are constants $k_i$ such that

$$h_1'''(t) \sim k_1 t^{\alpha-3} \text{ and } h_2'''(t) \sim k_2 t^{\alpha-3} \text{ as } t \to 0. \tag{10}$$

In principle, one can then read off the asymptotic behaviour of $I_2(x)$ as $x \to \infty$. However, although there are known results linking the asymptotic behaviour of Fourier transforms as $x \to 0$ to the behaviour of the integrand at $\infty$, (see e.g. [6]), there doesn't seem to be the necessary Abelian theorem in our context. We therefore formulate such a result, without striving for the greatest generality, in Proposition 2 below. The claimed result then follows easily.

2. Proofs

**Lemma 2.1:** $I_1(x) = O(x^{-3})$ as $x \to \infty$.

**Proof.** *Noting that, since $\alpha \in (1,2)$, $a := -\cos \alpha\pi/2$ and $b = \sin \alpha\pi/2$ are positive and $\beta := 1 - \alpha^{-1} \in$*

$(0, 1/2)$ we can make the change of variable $y = t^\alpha z$, to see that (7) holds with

$$g_1(t) = \int_0^{t^\alpha} e^{-ay}(t^\alpha - y)^{-\beta} \sin(by) dy$$

$$= t \int_0^1 e^{-azt^\alpha}(1-z)^{-\beta} \sin(bzt^\alpha) dz := tf_1(t), \text{ and}$$

$$g_2(t) = \int_0^{t^\alpha} e^{-ay}(t^\alpha - y)^{-\beta} \cos(by) dy$$

$$= t \int_0^1 e^{-azt^\alpha}(1-z)^{-\beta} \cos(bzt^\alpha) dz := tf_2(t).$$

Then, since it is easy to check that $tf_1(t)$ and its first two derivatives vanish both at 0 and $\infty$, we can integrate by parts 3 times to get

$$x^3 \int_0^\infty g_1(t) \cos(tx) dt = \int_0^\infty g_1'''(t) \sin(tx) dt = \int_0^\infty [tf_1'''(t) + 3f_1''(t)] \sin(tx) dt.$$

Then, making the change of variable $t = y/z^{\frac{1}{\alpha}}$, we can see that

$$\int_0^\infty |tf_1'''(t)| dt \leq k \int_0^\infty dt \int_0^1 e^{-azt^\alpha} \{zt^{\alpha-2} + z^2 t^{2\alpha-2} + z^3 t^{3\alpha-2}\}(1-z)^{-\beta} dz$$

$$= k \int_0^\infty e^{-ay^\alpha} dy \{y^{\alpha-2} + y^{2\alpha-2} + y^{3\alpha-2}\} \int_0^1 z^{1/\alpha}(1-z)^{-\beta} dz < \infty.$$

Similarly

$$\int_0^\infty |f_1''(t)| dt \leq k \int_0^\infty dt \int_0^1 e^{-azt^\alpha} \{zt^{\alpha-2} + z^2 t^{2\alpha-2}\}(1-z)^{-\beta} dz < \infty.$$

Virtually the same argument applies to $f_2$. (Note that although $(tf_2(t))' \to f_2(0) \neq 0$ as $t \to 0$, we only need $(tf_2(t))' \sin tx \to 0$ to justify the second integration by parts, and this does hold.) ∎

Turning to $I_2(x)$, we note that (9) holds with $h_1(t) = e^{-at^\alpha} \cos bt^\alpha$ and $h_2(t) = -e^{-at^\alpha} \sin bt^\alpha$. Moreover it is clear that (10) holds with $k_1 = a\alpha(\alpha-1)(2-\alpha)$ and $k_2 = b\alpha(\alpha-1)(2-\alpha)$. So we need

**Proposition 2.2:** Suppose $h \in C^3(0, \infty)$, $\sup_{0<t<\infty} |t^{3-\alpha} h'''(t)| < \infty$, and

$$h'''(t) \backsim t^{\alpha-3} \text{ as } t \to 0, \tag{11}$$

where $1 < \alpha < 2$. Assume also that each of $\lim_{t\to\infty} h(t), \lim_{t\to\infty} h'(t)$ and $\lim_{t\to\infty} h''(t)$ are 0. Then (i) if $\lim_{t\to 0} th(t) = \lim_{t\to 0} h'(t) = 0$,

$$\int_0^\infty h(t) \cos(tx) dt \backsim l_1 x^{-\alpha-1} \text{ as } x \to \infty,$$

where $l_1 = \int_0^\infty y^{\alpha-3} \sin y \, dy = \frac{-\pi}{2\Gamma(3-\alpha) \cos \alpha\pi/2}$, and
(ii) if $\lim_{t\to 0} h(t) = \lim_{t\to 0} th'(t) = 0$,

$$\int_0^\infty h(t) \sin(tx) dt \backsim l_2 x^{-\alpha-1} \text{ as } x \to \infty,$$

where $l_2 = (2-\alpha)^{-1}\lim_{x\to\infty}\int_0^x y^{\alpha-2}\sin y\,dy = \frac{\pi}{2\Gamma(3-\alpha)\sin\alpha\pi/2}$.

**Proof.** (i) First note that we can integrate (11) to see that $\lim_{t\to 0} th''(t) = 0$. Writing $h'''(t) = t^{\alpha-3}j(t)$ this, together with the other assumptions, enables us to integrate by parts three times to get

$$x^{\alpha+1}\int_0^\infty h(t)\cos tx\,dt = x^{\alpha-2}\int_0^\infty h'''(t)\sin tx\,dt$$

$$= \int_0^\infty y^{\alpha-3}j(y/x)\sin y\,dy \to l_1,$$

where we have used dominated convergence. The expression for $l_1$ comes from formula 3.761(4), p 420 of [5].

(ii) Two integrations by part give

$$x^{\alpha+1}\int_0^\infty h(t)\sin tx\,dt = -x^{\alpha-1}\int_0^\infty h''(t)\sin tx\,dt = -x^{\alpha-1}(I_1(x)+I_2(x)),\text{ where}$$

$$-x^{\alpha-1}I_1(x) = -x^{\alpha-1}\int_{\pi/2x}^\infty h''(t)\sin tx\,dt$$

$$= x^{\alpha-2}\int_{\pi/2x}^\infty h'''(t)\cos tx\,dt \to \int_{\pi/2}^\infty y^{\alpha-3}\cos y\,dy,$$

by the same argument. However we also have $h''(t) \sim (\alpha-2)^{-1}t^{\alpha-2}$ so we get

$$-x^{\alpha-1}I_2(x) = -x^{\alpha-1}\int_0^{\pi/2x} h''(t)\sin tx\,dt$$

$$= -\int_0^{\pi/2}(y/x)^{2-\alpha}h''(y/x)y^{\alpha-2}\sin y\,dy$$

$$\to (2-\alpha)^{-1}\int_0^{\pi/2} y^{\alpha-2}\sin y\,dy.$$

Again the expression for $l_2$ follows from 3.761(4), p 420 of [5]. ∎

**Theorem 2.3:** *(3) holds.*

**Proof.** Since both of $k_i^{-1}h_i(t)$ satisfy the assumptions of Proposition 1, the result follows once we check that

$$\frac{k_1 l_1 + k_2 l_2}{\pi} = \frac{\alpha(\alpha-1)(2-\alpha)}{\pi}(-l_1\cos\alpha\pi/2 + l_2\sin\alpha\pi/2)$$

$$= \frac{\alpha(\alpha-1)(2-\alpha)}{\pi}\cdot\frac{\pi}{\Gamma(3-\alpha)} = \frac{1}{\Gamma(-\alpha)}.$$

∎

**Remark 1:** Although the asymptotic behaviour of $P(S_1 > x)$ is known, and of course follows from our result, we point out it also follows easily from (4). Indeed, writing $1/\alpha = \eta$ and making an obvious change of variable gives

$$\int_0^\infty e^{-\lambda x}s(x)dx = e^{\lambda^\alpha}(1 - \frac{1}{\Gamma(\eta)}\int_0^{\lambda^\alpha} y^{\eta-1}e^{-y}dy).$$

It is obvious that $\int_0^{\lambda^\alpha} y^{\eta-1} e^{-y} dy = \alpha\lambda + O(\lambda^2)$ and hence

$$\int_0^\infty e^{-\lambda x} s(x) dx = 1 - \frac{\alpha\lambda}{\Gamma(\alpha)} + \lambda^\alpha + O(\lambda^2) \text{ as } \lambda \downarrow 0,$$

and then (2) follows from Theorem 8.1.6 in [3].

**Remark 2:** The asymptotic behaviour of $P(S_1 > x)$ is the same, modulo the value of the constant, for any stable process of index $\alpha \in (1,2)$ which has positive jumps. It is therefore possible that the result of our Theorem also holds in this generality. However, in the absence of any specific information about the Wiener-Hopf factors, it is not clear how this could be established.